\documentclass[pmlr]{jmlr}

\usepackage{longtable}
\usepackage{bbm}

\usepackage{booktabs}
\usepackage[load-configurations=version-1]{siunitx} 

\theorembodyfont{\upshape}
\theoremheaderfont{\scshape}
\theorempostheader{:}
\theoremsep{\newline}

\jmlryear{2027}
\jmlrworkshop{Learning for Dynamics and Control}

\title[Learning-Based Hierarchical Planning and Tracking]{Certified Hierarchical Planning and Tracking: A Learning-Based Approach}

\author{\Name{Weijian Li\nametag{$^{1}$}} \Email{wli26@nd.edu}\\
\Name{Vince Kurtz\nametag{$^{2}$}}        \Email{vkurtz1@depaul.edu} \\
\Name{Panos J. Antsaklis\nametag{$^{1}$}} \Email{antsaklis.1@nd.edu} \\
\Name{Hai Lin\nametag{$^{1}$}}            \Email{hlin1@nd.edu}\\
\addr
$^{1}$Department of Electrical Engineering, University of Notre Dame\\
$^{2}$School of Computing, DePaul University
}

\begin{document}

\maketitle

\begin{abstract}

In this paper, we develop a hierarchical control framework for nonlinear dynamical systems with safety constraints. Specifically, a low-fidelity model is introduced for trajectory planning, while a tracking controller is designed for the high-fidelity system. We first characterize the tracking error bound (TEB) using a storage function. Then, we formulate the synthesis of the tracking controller, the storage function, and the TEB as an optimization problem. Following that, we develop an offline learning approach to approximately solve the optimization problem, by incorporating behavior cloning and satisfiability modulo theories (SMT) solver.
We further integrate the learned controller and the TEB into online tracking of the high-fidelity system.
This framework enables the control of high-dimensional nonlinear systems with safety guarantees, while maintaining computational efficiency for real-time implementation.
Finally, experiments are provided for illustration.
\end{abstract}


\section{Introduction}
\label{sec:intro}

Designing controllers with performance and safety guarantees is of fundamental importance for modern engineering systems, such as autonomous vehicles and robotics.
However, controller synthesis and real-time implementation remain challenging due to the nonlinear and high-dimensional dynamics of these systems \citep{hsu2023safety, cohen2024safety}.
Hierarchical control provides a promising paradigm to address these challenges through a layered architecture for planning and tracking \citep{girard2009hierarchical, schweidel2022safe, singh2018robust, kousik2020bridging, chen2021fastrack}.
In particular, a lower-fidelity model (planner) is introduced for online planning, while the high-fidelity model (tracker) follows the planned trajectory via a tracking controller.
Such a separation provides flexibility and computational tractability for the real-time control of complex dynamical systems \citep{herbert2017fastrack, benders2025embedded}.

In practice, it is significant to provide formal guarantees for hierarchical control, as the mismatch between low- and high-fidelity models generally results in tracking errors.
In \cite{girard2009hierarchical, niu2026approximate}, tracking errors were characterized through approximate simulation, where simulation functions were employed to establish TEBs between the concrete systems and their abstractions.
Alternatively, Hamilton-Jacobi reachability analysis was employed in \cite{herbert2017fastrack, chen2021fastrack, jeong2024parameterized} to compute TEBs and synthesize safe tracking controllers.
In \cite{smith2019continuous, majumdar2017funnel, singh2018robust, schweidel2022safe}, storage functions have been employed to characterize invariant tracking error sets, with the synthesis of these functions and tracking controllers formulated as sum-of-squares (SOS) programs. However, these SOS-based approaches are typically limited to polynomial dynamical systems and become computationally expensive for large-scale systems.

Recent advances in learning-based control have provided flexible approaches for stability analysis and controller synthesis of complex nonlinear systems. For instance, neural Lyapunov functions have been developed to estimate regions of attraction of nonlinear systems \citep{liu2025physics, chang2019neural, zhang2023compositional}. Neural networks have also been employed to approximate optimal control policies for nonlinear systems \cite{meng2024physics, lu2024mpc}. In \cite{dawson2023safe, sun2021learning, ma2022joint},
controllers and certificate functions have been jointly learned to simultaneously achieve control objectives and provide stability or safety guarantees.

\noindent\textbf{Our Contributions.} In this paper, we develop a learning-based framework for hierarchical planning and tracking of nonlinear dynamical systems with safety constraints. Our main contributions are at least threefold. 
(i) By characterizing the tracking error through a forward invariant sublevel set of an appropriate storage function, we formulate the synthesis of the tracking controller, the storage function, and the TEB as an optimization problem, thereby providing a systematic approach to the hierarchical design.
(ii) We propose a learning-based approach for the optimization problem by integrating behavior cloning with SMT-based falsification. The method is applicable to general nonlinear dynamical systems, and explicitly incorporates certifications to provide performance guarantees.
(iii) We develop an online planner-tracker architecture that integrates the learned tracking controller and TEB to ensure safe tracking of the planned trajectory.
We validate the effectiveness of the proposed framework through examples from real-world applications.

The remainder of this paper is organized as follows. Section~\ref{sec:formulation} formulates the hierarchical planning and tracking problem. Section~\ref{sec:frame} develops the theoretical framework.
Section~\ref{sec:offline} presents the learning-based approach, while Section~\ref{sec:online} focuses on the online implementation. Finally, Section~\ref{sec:experiment} presents numerical experiments, and Section~\ref{sec:conclusion} concludes the paper.

\section{Problem Setup}
\label{sec:formulation}

Consider a high-fidelity model, referred to as the tracking model, described by a control-affine system of the form
\begin{equation}
\label{dyn:high}
\Sigma:~
\dot x = f(x, d) + g(x, d)u,
\end{equation}
where $x \in \mathcal X \subset \mathbb R^n$ is the system state,  
$u \in  \mathcal U \subset \mathbb R^m$ is the control input, 
$d \in  \mathcal D \subset \mathbb R^p$ is the disturbance. Both $f: \mathbb R^n \times \mathbb R^p \rightarrow \mathbb R^n$ and $g: \mathbb R^n \times \mathbb R^p \rightarrow \mathbb R^{n \times m}$ are smooth, $\mathcal X$ and $\mathcal U$ are constraint sets imposed on the state and input, and $\mathcal D$ is compact.

We refer to $\Sigma$ as the concrete system, which may be nonlinear and high-dimensional. To control such a complex system, we introduce a low-fidelity model  within a hierarchical control framework \citep{schweidel2022safe, herbert2017fastrack, singh2018robust, siriya2025towards},  which serves as an abstraction of $\Sigma$ and is given by
\begin{equation}
\label{dyn:low}
\hat \Sigma:~
\dot z = \hat f(z, v),
\end{equation}
where  $z \in \mathcal Z \subset \mathbb R^{\hat n}$,
$v \in \mathcal V \subset \mathbb R^{\hat m}$ with $\mathcal V$ being compact, and
$\hat f: \mathbb R^{\hat n} \times \mathbb R^{\hat m} \rightarrow \mathbb R^{\hat n}$.
In general, $\hat \Sigma$ has a lower dimension than $\Sigma$ and is used as a planning model for online trajectory generation.
The kinematics or simplified dynamics of $\Sigma$ may be employed as the planner.
A controller $u$ is then synthesized offline for \eqref{dyn:high} to track the trajectory generated by $\hat \Sigma$. 
To characterize the tracking performance, we define the error variable as
\begin{equation}
\label{error}
e(t) := \phi(x,z)(x - \pi(z)),
\end{equation}
where  $\pi: \mathbb R^{\hat n} \rightarrow \mathbb R^n$ is a $\mathcal C^1$ (continuously differentiable) map that lifts the planning state \(z\) to the state space of $\Sigma$, and $\phi: \mathbb R^n \times \mathbb R^{\hat n} \rightarrow \mathbb R^{\tilde n \times n}$ is $\mathcal C^1$, providing additional flexibility.

It follows from \eqref{dyn:high} and \eqref{dyn:low} that
\begin{equation}
\begin{aligned}
\label{dyn:error}
\dot e := \tilde f(e, u, v, d) 
= f_e\bigl(e, v, d\bigr) + g_e\big(e, d\big) u,
\end{aligned}
\end{equation}
where 
\begin{equation*}
\begin{aligned}
f_e(e, v, d) =&~ \Big[
\frac{\partial \phi}{\partial x} f(x, d) + \frac{\partial \phi}{\partial z} \hat f(z, v)\Big]\big(x - \pi(z)\big) 
+ \phi(x, z)\Big[f(x, d) - \frac{\partial \pi}{\partial z}  
\hat f(z, v)\Big], \\
g_e(e, d) =&~\frac{\partial}{\partial x}\big[\phi(x, z) \big(x- \pi(z)\big)\big]
g(x, d).
\end{aligned}
\end{equation*}
Throughout this paper, we assume that $\tilde f$ admit a representation of the form  \eqref{dyn:error},
depending only on $e$, $u$, $v$ and $d$. 
This is consistent with the formulations in \cite{singh2018robust, schweidel2022safe, siriya2025towards}, and many practical applications admit error dynamics of this form.
In Appendix \ref{apd:first}, we provide two examples for illustration.

Based on the above framework, the goal of this paper is threefold: 
(i) design a tracking controller for \eqref{dyn:high}, 
(ii) characterize the tracking error between $\Sigma$ and $\hat \Sigma$,
and (iii) incorporate the resulting controller and the error characterization into the planner-tracker architecture for online trajectory planning and tracking.

\section{Hierarchical Planning and Tracking}
\label{sec:frame}

In this section, we first characterize the TEB using a storage function. We then formulate the synthesis of the tracking controller, TEB and storage function as an optimization problem. Finally, we present a learning-based framework for hierarchical planning and tracking.

\subsection{Tracking Error Bound}

Consider a control law for \eqref{dyn:high} of the form
\begin{equation}
\label{contrl:high}
u = \kappa(e, v),
\end{equation}
where $\kappa: \mathbb R^{\tilde n} \times \mathcal V \rightarrow \mathcal U$.
Plugging \eqref{contrl:high} into \eqref{dyn:error}, we obtain
\begin{equation}
\begin{aligned}
\label{dyn:error:ctrl}
\dot e = f_e\bigl(e, v, d\bigr) + g_e\big(e, d\big) \kappa(e, v).
\end{aligned}
\end{equation}

For dynamics \eqref{dyn:error:ctrl}, a set $\mathcal O$ is a TEB from an initial set $\mathcal I$ if, for any $e(0) \in \mathcal I$, $v \in \mathcal V$ and 
$d \in \mathcal D$, it holds that $e(t) \in \mathcal O$ for all $t \ge 0$.
The following proposition, proved in Appendix \ref{apd:second}, provides a sufficient condition for establishing a TEB.
\begin{proposition}
\label{prop:Lyap}
Consider the error dynamics \eqref{dyn:error:ctrl}. 
Suppose that there exist a $\mathcal C^1$ function $V: \mathbb R^{\tilde n} \rightarrow \mathbb R_+$, $\kappa: \mathbb R^{\tilde n} \times \mathcal V \rightarrow \mathcal U$, and a constant $\gamma > 0$ such that
\begin{equation}
\begin{aligned}
\label{lyap:cond}
V(0) &= 0, ~V(e) > 0, ~\forall e \neq 0, ~\text{and}~
\lim_{\|e\| \to \infty}V(e) = \infty,\\
V(e) \ge \gamma &\Longrightarrow 
\frac{\partial V}{\partial e} \cdot \Bigl(f_e(e,v,d) + g_e(e,d)\kappa(e,v)\Bigr) < 0, 
~\forall e \in \mathbb R^{\tilde n}, \forall v\in \mathcal V, \forall d\in \mathcal D.
\end{aligned}
\end{equation}
For the sublevel set $\Omega^V_{\gamma} := \{e \in \mathbb R^{\tilde n}: V(e) \le \gamma\}$, the following statements hold.
\begin{enumerate}
\renewcommand{\labelenumi}{(\alph{enumi})}
\item The set $\Omega^V_{\gamma}$ is forward invariant, i.e., $e(t) \in \Omega^V_{\gamma}$ if $e(0) \in \Omega^V_{\gamma}$.
\item For all $e(0) \in \mathbb R^{\tilde n}$,  $e(t)$ remains bounded, and moreover, 
$\lim_{t \to \infty}{\rm dist}(e(t), \Omega_\gamma^V) = 0$.
\end{enumerate}
\end{proposition}

\begin{remark}
Proposition \ref{prop:Lyap} ensures that $\Omega_\gamma^V$ is a TEB and that the tracking error converges toward this set. 
In contrast to \cite{smith2019continuous}, we explicit account for disturbances in the high-fidelity model \eqref{dyn:high} and the resulting error dynamics \eqref{dyn:low}. 
Compared with \cite{schweidel2022safe}, where $\dot V < 0$ is imposed only on the boundary $V(e) = \gamma$, Proposition \ref{prop:Lyap} requires the condition over the region $V(e)\ge\gamma$. This facilitates the subsequent learning procedure, since the decrease condition allows sampling outside $\Omega^V_{\gamma}$ instead of sampling on its boundary.
\end{remark}

\subsection{Optimization-Based Formulation}

To achieve a conservative TEB, we consider a optimization problem with the volume of $\Omega^V_{\gamma}$ as the cost function. 
Then, the problem is formulated as
\begin{equation}
\begin{aligned}
\label{form}
\min_{V, \kappa}~\text{volume}(\Omega^V_{\gamma}),~~
{\rm s.t.}~\eqref{dyn:error:ctrl}-\eqref{lyap:cond}.
\end{aligned}
\end{equation}

\begin{remark}
We should mention that the set $\Omega^V_\gamma$ is equivalent, up to a scaling factor, to the set $\{e: \alpha V(e) \le \alpha \gamma\}$, where $\alpha > 0$.
To eliminate the redundancy, without loss of generality, we set $\gamma = 1$ throughout the rest of this paper \citep{singh2018robust}.
\end{remark}

In practice, it is difficult or even impossible to compute $\text{volume}(\Omega^V_{\gamma})$.
Given $\gamma = 1$, we approximate the set $\Omega^V_{\gamma}$ by an ellipsoid 
$\mathcal E := \{e: e^\top E e \le 1\}$ \citep{majumdar2017funnel},
where $E \in \mathbb R^{\tilde n \times \tilde n}$ is positive definite.
We require that $e(t) \in \mathcal E$ if $e(t) \in \Omega^V_{\gamma}$.
Then $\mathcal E$ characterizes the TEB as that of $\Omega^V_\gamma$.
Note that the volume of $\mathcal E$ is proportional to $1/{\sqrt{\det(E)}}$.
Thus, problem \eqref{form} can be reformulated as
\begin{equation}
\begin{aligned}
\label{reform}
\min_{V, \kappa, E}&~ -\log\det(E) \\
\text{s.t.}&~V(0) = 0, 
~V(e) > 0, ~\forall e \neq 0, 
~\lim_{\|e\| \to \infty}V(e) = \infty,\\
&~V(e) \ge 1 \Longrightarrow \frac{\partial V}{\partial e}
\cdot \Bigl(f_e(e,v,d) +
g_e(e,d)\kappa(e,v)\Bigr) < 0, \\
&~V(e) \le 1 \Longrightarrow e^\top E e \le 1.
\end{aligned}
\end{equation}

\begin{remark}
By solving problem \eqref{reform}, we jointly seek a tracking controller $\kappa$, a storage function $V$ and an ellipsoidal TEB $\mathcal E$. 
In \cite{smith2019continuous, singh2018robust, majumdar2017funnel}, SOS programs were employed to synthesize tracking controllers together with TEB or funnels, but relied on polynomial representations of the system dynamics. In contrast, we develop a learning-based approach that can handle more general dynamics.
\end{remark}

\subsection{Proposed Framework}

We are now ready to present the proposed framework for approximately solving \eqref{reform} offline and integrating the resulting tracking controller and TEB into the planner-tracker architecture for online implementation.
Problem \eqref{reform} provides a unified formulation for addressing the objectives stated in Section~\ref{sec:formulation}.
In general, however, it is difficult to be solved for the class of nonlinear control-affine systems considered in \eqref{dyn:high}.
To address this challenge, we parameterize the tracking controller $\kappa$ and storage function $V$ using neural networks and seek an approximate solution to \eqref{reform} through an offline learning-based approach.

\begin{figure}[t]
\centering
\includegraphics[
width=\textwidth,
trim=3.2em 19.8em 18em 5em,
clip
]{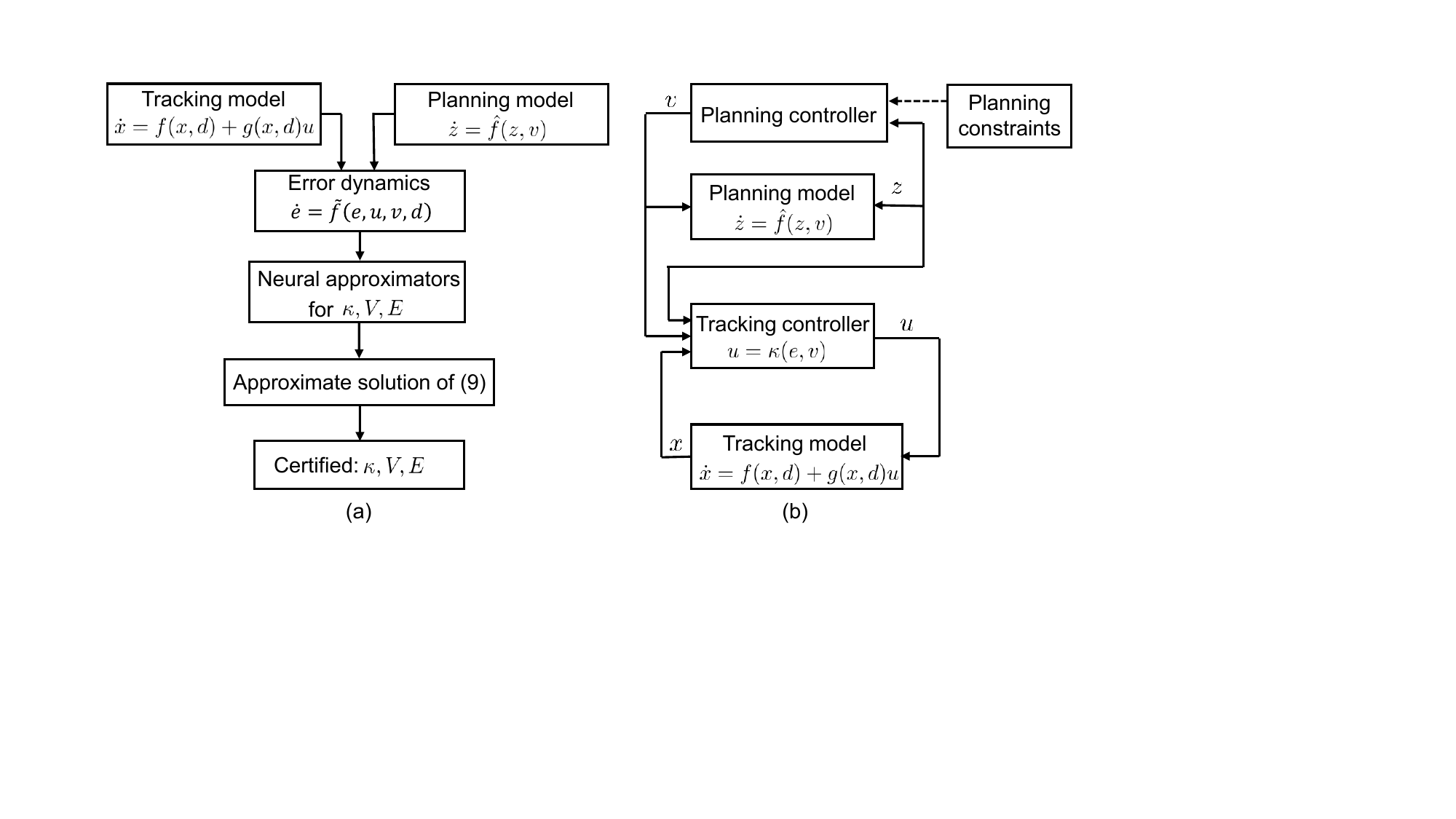}
\caption{Overview of our framework: 
(a) Offline synthesis;
(b) Online implementation.}
\label{Fig1}
\end{figure}

For online implementation, given the TEB $\mathcal E$ and the state constraint set $\mathcal X$ for \eqref{dyn:high}, we first compute an admissible state set $\mathcal Z_{\rm safe}$  for the planning model \eqref{dyn:low}, such that any trajectory generated within $\mathcal Z_{\rm safe}$ can be safely tracked by \eqref{dyn:high}.
Using the planning model \eqref{dyn:low}, the admissible state set $\mathcal Z_{\rm safe}$, and the input constraint $\mathcal V$, the control input $v(t)$ for \eqref{dyn:low} can be computed.  
Following that, the learned control policy $u$ is applied to \eqref{dyn:high}.
Fig.~\ref{Fig1} demonstrates the offline synthesis and online implementation of the proposed planner-tracker architecture.

\section{Offline Learning}
\label{sec:offline}

\subsection{Parameterization and Loss Function}

To approximately solve problem \eqref{reform}, we begin with parameterizing the tracking controller as $u = \kappa_{\theta_u}(e, v)$, where $\kappa_{\theta_u}(\cdot)$ is implemented by a feedforward neural network with parameters $\theta_u$.
In practice, a saturation layer can be applied to the output of the neural network to enforce hypercube constraints on $u$. To ensure positive definiteness and radial unboundedness, we parameterize the storage function $V$ as
\begin{equation*}
V_{\theta_V}(e) := C^\top(e; \theta_V) C(e; \theta_V) + \delta_V e^\top e,
\end{equation*}
where $C(\cdot; \theta_V)$ is a neural network, $\theta_V$ denotes the trainable parameters, and $\delta_V > 0$ is a user-defined hyperparameter. For the ellipsoidal TEB $\mathcal E$, we parameterize $E$ as $E= S^\top S + \delta_E I_{\tilde n}$,  where $S \in \mathbb R^{\tilde n \times \tilde n}$ is the decision variable, and $\delta_E > 0$ is a user-defined hyperparameter.

Consider a training dataset $\mathcal D:=\{(e_i,v_i,d_i)\}_{i=1}^{N}$, where $e_i$, $v_i$, and $d_i$ are sampled from their training domains. We encode the certificate conditions in \eqref{reform} as empirical loss functions using suitable penalty terms.
In particular, the resulting loss function consists of three components: a storage function loss, an  ellipsoidal bound loss, and a behavior cloning loss, as follows.

To enforce $\dot V < 0$ for $V(e) \ge 1$ and $V(0) = 0$, the empirical loss for $V_{\theta_V}$ is defined as
\begin{equation}
\label{loss:V}
\mathcal L_V := \frac 1N \sum_{i = 1}^N \Big(\mathbbm{1}_{V_{\theta_V}(e_i) \ge 1} \max\big(0, \dot V_{\theta_V}(e_i, v_i, d_i) + \epsilon\big) \Big) + \alpha V^2_{\theta_V}(0),
\end{equation}
where  $\alpha > 0$ is a weight parameter, $\dot V_{\theta_V}$ denotes the time derivative of $V_{\theta_V}$ along the error dynamics \eqref{dyn:error:ctrl}, $\epsilon > 0$ is a small margin introduced to encourage $\dot V_{\theta_V} < 0$, and the indicator function $\mathbbm{1}_{V_{\theta_V}(e_i) \ge 1}(\cdot)$ ensures that the penalty is activated only if $V_{\theta_V}(e_i) \ge 1$.

To promote $\Omega_1^V\subseteq\mathcal E$ and reduce the size of $\mathcal E$, we define the ellipsoidal bound loss as
\begin{equation}
\label{loss:E}
\mathcal L_E := \frac 1N \sum_{i = 1}^N \Big(\mathbbm{1}_{V_{\theta_V}(e_i)\le 1}
\max \big(0,e_i^\top E e_i - 1 \big) \Big) - \beta \log \det(E),	
\end{equation}
where $\beta > 0$ is a weight parameter and the indicator function restricts the penalty to $\Omega_1^V$.

The losses $\mathcal L_V$ and $\mathcal L_E$ account for certificate conditions in \eqref{reform}.
However, the performance of the learned controller may be poor, especially during the early stages of training, as unreliable control actions may be generated.
To address this issue, we introduce a behavior cloning loss that encourages the learned controller to imitate expert policies. 
Consider a set of expert demonstrations $\mathcal D_{BC} := \{(e_j, v_j, u^{\rm ep}_j)\}_{j = 1}^{N_{\rm BC}}$,
where $u_j^{\rm ep}$ denotes the expert control input corresponding to $(e_j, v_j)$.
The behavior cloning loss is then defined as
\begin{equation}
\label{loss:BC}
\mathcal L_{\text{BC}} = \frac {1}{N_{\rm BC}} \sum_{j = 1}^{N_{\rm BC}} \|\kappa_{\theta_u}(e_j, v_j) - u^{\text{ep}}_j\|^2.
\end{equation}
We note that $u^{\text{ep}}_j$ is derived by solving a finite-horizon model predictive control (MPC) problem based on  \eqref{dyn:error}. For the MPC implementation, we linearize and discretize the continuous-time nonlinear error dynamics.

Combining \eqref{loss:V}, \eqref{loss:E} and \eqref{loss:BC}, the overall loss function is defined as
\begin{equation}
\label{loss}
\mathcal L(\theta_u, \theta_V, S) = \omega_V \mathcal L_V + \omega_E \mathcal L_E + \omega_{BC}\mathcal L_{\text{BC}},
\end{equation}
where $\omega_V$, $\omega_E$ and $\omega_{\text{BC}}$ are user-defined positive weights.

\begin{remark}
By minimizing the loss \eqref{loss}, the tracking controller, storage function, and TEB can be jointly learned. The loss $\mathcal L_V$ in \eqref{loss:V} is similar to the neural Lyapunov loss in \cite{chang2019neural,zhou2022neural,gaby2022lyapunov}, but differs in that the penalty is imposed only when
$V(e) \ge 1$. 
The loss $\mathcal L_E$ provides a tractable way to integrate the condition $\Omega_1^V\subseteq\mathcal E$ with the TEB minimization. Behavior cloning, explored in \cite{sun2021learning,ahn2023model,dawson2023safe}, is incorporated to facilitate the training of the controller.
\end{remark}

\subsection{Training and Verification}

To verify the certificate conditions beyond the training samples, we use an SMT solver \citep{gao2013dreal, chang2019neural, zou2025analytical} as a falsifier to search for counterexamples to the learned solution $(\kappa_{\theta_u}, V_{\theta_V}, E)$. To be specific, 
we formulate the counterexample conditions as
\begin{equation}
\label{SMT:falsifier}
\Phi(e,v,d):=\Phi_V(e,v,d)\vee \Phi_E(e),
\end{equation}
where 
\begin{equation*}
\begin{aligned}
\Phi_V(e, v, d) &:= \big(V_{\theta_V}(e)\ge 1\big)
\wedge \big(v\in\mathcal V\big)
\wedge \big(d\in\mathcal D\big) \wedge (\dot V_{\theta_V} \ge 0), \\
\Phi_E(e) &:= \big(V_{\theta_V}(e) \le 1\big)
\wedge \big(e^\top E e > 1\big).
\end{aligned}
\end{equation*}

If $\Phi$ is satisfiable, the corresponding $(e,v,d)$ provides a counterexample, which is added to the training set. 
Otherwise, no counterexample exists within the predefined search domain, up to the numerical precision of the SMT solver \cite{gao2012delta}.

The training and falsification steps are performed iteratively in a counterexample-guided manner.
We summarize the offline synthesis procedure in Algorithm \ref{alg:offline}.

{\SetAlCapHSkip{0pt}
\begin{algorithm}[t]
\hrule height 1.0pt
\vspace{0.9ex}
\caption{\noindent Joint Learning of the Tracking Controller, Storage Function, and TEB}
\label{alg:offline}
\vspace{0.4ex}
\hrule
\vspace{0.6ex}

\KwIn{
Tracking model $\Sigma$, planning model $\hat{\Sigma}$, error mapping \eqref{error}, maximum number of iterations $N_{\max}$, and initial parameters $\theta_u$, $\theta_V$, and $S$;}

\KwOut{
Controller $\kappa_{\theta_u}(e,v)$, storage function $V_{\theta_V}(e)$ and TEB $\mathcal E$;}

Construct the error dynamics \eqref{dyn:error:ctrl};

Generate certificate samples $\mathcal D$ and expert demonstrations $\mathcal D_{\rm BC}$;

\For{$k=1,\ldots,N_{\max}$}{

Train $(\theta_u,\theta_V,S)$ by minimizing \eqref{loss} using
$\mathcal D$ and $\mathcal D_{BC}$;

Obtain the candidate $(\kappa_{\theta_u},V_{\theta_V},E)$;

Run the SMT falsifier using \eqref{SMT:falsifier};

\eIf{UNSAT}{\Return
	$(\kappa_{\theta_u},V_{\theta_V},E)$;}{
	Collect counterexamples $\mathcal C$;
	
	Update
	$\mathcal D = \mathcal D \cup\mathcal C$;
}}
\Return
$(\kappa_{\theta_u},V_{\theta_V}, \mathcal E)$;
\vspace{0.9ex}
\hrule height 1.0pt
\end{algorithm}}

\section{Online Implementation}
\label{sec:online}

In this section, we incorporate the learned controller $\kappa_{\theta_u}$ and TEB $\mathcal E$ into the planner-tracker architecture for online implementation.
We use $\mathcal E$ to tighten the state constraints of \eqref{dyn:low}, such that the trajectory generated by $\hat \Sigma$ can be safely tracked by \eqref{dyn:high}.

Suppose that $e(0) \in \Omega_1^V$. By Proposition \ref{prop:Lyap}, $e(t) \in \mathcal E$ for all $t \ge 0$.
Given the planning state $z$, we define the set of tracking states associated with $z$ and $\mathcal E$ as
\begin{equation*}
\mathcal R(z) :=
\big\{x\in\mathbb R^n:
\phi(x,z)\bigl(x-\pi(z)\bigr)\in\mathcal E
\Big\}.
\end{equation*}
To ensure $x \in \mathcal X$, we require $\mathcal R(z) \subset \mathcal X$. Accordingly, the admissible state set for \eqref{dyn:low} is defined as
\begin{equation*}
\label{planner:constraint}
\mathcal Z_{\text{safe}} := \big\{z\in\mathbb R^{\hat n} : \mathcal R(z)\subseteq\mathcal X \big\}.
\end{equation*}
When $\phi(x,z)= I$, $\mathcal Z_{\text{safe}} \!=\! \{z \in \mathbb R^{\hat n}: \pi(z)\in \mathcal X \ominus\mathcal E\}$, where $\ominus$ is the Pontryagin difference. 

For online planning, we linearize and discretize the planning model \eqref{dyn:low}, and employ MPC subject to the admissible state constraint $z(t)\in\mathcal Z_{\mathrm{safe}}$ to compute the planning input $v(t)$. The learned controller $u(t)=\kappa(e(t),v(t))$ is then applied to the tracking model \eqref{dyn:high}.

\section{Experiments}
\label{sec:experiment}

In this section, we evaluate the performance of the proposed approach through three examples from real-world applications.
In particular, we consider hierarchical regulation problems, in which \eqref{dyn:low} plans a trajectory toward a fixed goal, while  \eqref{dyn:high} follows the planned trajectory. Given \eqref{error}, we first derive the error dynamics \eqref{dyn:error}.
For numerical implementation, we discretize the continuous-time dynamics \eqref{dyn:high}, \eqref{dyn:low}, and \eqref{dyn:error} with a sampling period $\Delta T$. To generate the training data, we randomly sample a finite number of regulation tasks. For each task, the initial and goal states of the planner \eqref{dyn:low} are sampled, and 
the inputs are computed by solving a finite-horizon MPC problem, yielding a sequence of state-input pairs $\{(z_k,v_k)\}$.
At each planner state $z_k$, a tracking error $e_k$ is independently sampled. For each tuple $(e_k,z_k,v_k)$, an expert control input $u^{{\rm ep}}_k$ is obtained by solving a finite-horizon MPC problem based on the linearized and discretized tracking dynamics \eqref{dyn:high}. 
The resulting samples are then used to construct the training dataset.
Following that, we parameterize the tracking controller $u$, storage function $V$, and TEB $\mathcal E$, and jointly optimize the parameters by minimizing \eqref{loss} using stochastic gradient descent algorithm, as shown in Algorithm~1.
Finally, we incorporate the learned $u$ and $\mathcal E$ into the online planner--tracker architecture, as discussed in Section \ref{sec:online}.
The source code and parameter settings for all experiments are available at \url{https://github.com/WeijianLi-Maker/hierarchical-planning-tracking}.

\subsection{4D Car with a 2D Single-Integrator Planner}
\label{sec:example:4Dcar}

\begin{figure}[t]
\centering
\includegraphics[width=0.9\textwidth]{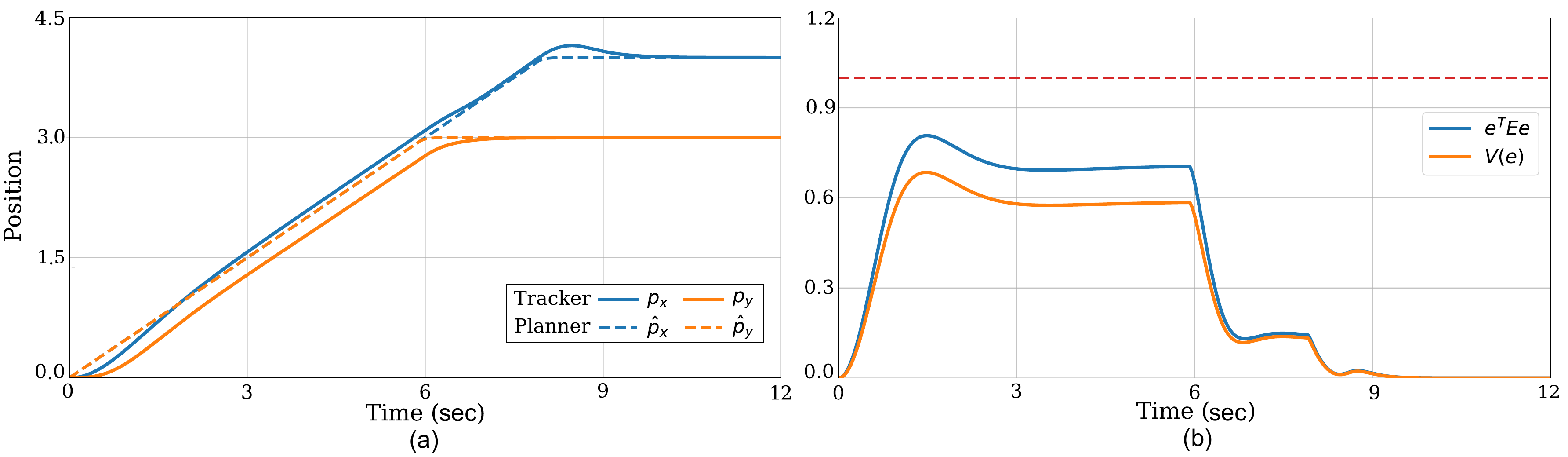}
\caption{(a) Trajectories of the planning and tracking models, where the input to \eqref{dyn:low} is computed by MPC and the input to \eqref{dyn:high} is generated by the learned controller. (b) Trajectories of $V(e)$ and  $e^\top E e$.}
\label{Fig2}
\end{figure}

Consider a car moving on a plane. The high-fidelity model is given by
$\dot p_x = \nu \cos(\psi)$, $\dot p_y = \nu \sin(\psi)$, $\dot \psi = \omega$ and $\dot \nu = a$, where $p_x$ and $p_y$ denote the position, $\psi$ is the heading angle, $\nu$ is the speed, $\omega$ and $a$ denote the angular velocity and acceleration inputs.
The planning model is $\dot{\hat p}_x = \hat \nu_x$ and $\dot{\hat p}_y = \hat \nu_y$, where $\hat \nu_x$ and $\hat \nu_y$ are the planning inputs.
Following \cite[Appendix A]{singh2018robust}, we define the relative error state as
$e = [(p_x-\hat p_x)\cos\psi+(p_y-\hat p_y)\sin\psi, -(p_x-\hat p_x)\sin\psi+(p_y-\hat p_y)\cos\psi, \nu]^\top$.
We impose box constraints on $\omega$ and $a$.
Fig.~\ref{Fig2} presents the tracking performance and the evolution of $V(e)$ and $e^\top E e$. 
The results demonstrate the effectiveness of the learned tracking controller. Moreover, the tracking error $e(t)$ remains within the sublevel set $\Omega^V_1$, which is contained in $\mathcal E$.
We provide additional simulation results under different goal states in Appendix~\ref{apd:forth} to further evaluate the tracking performance of the learned controller.

\subsection{Bicycle Model with a Dubin's Planner}

\begin{figure}[t]
\centering
\includegraphics[width=0.9\textwidth]{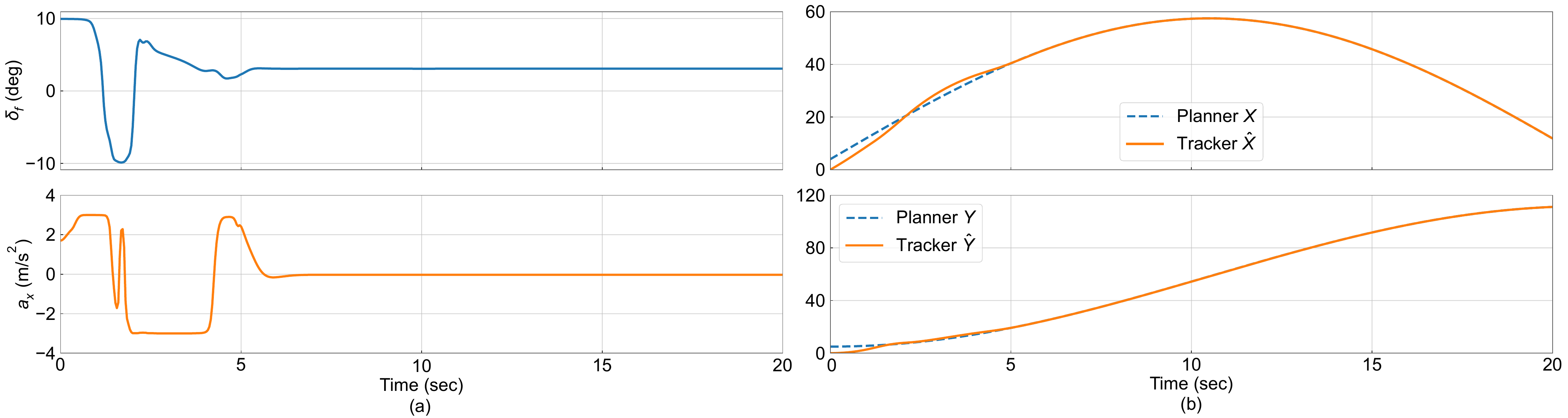}
\caption{(a) Control inputs of the tracking model. (b) Trajectories of the planning and tracking models with the tracking input generated by the learned controller. }
\label{Fig3}
\end{figure}

We revisit the bicycle tracking model and the Dubin's planner in Appendix \ref{append:A:2}. Similar to \cite{schweidel2022safe}, we set
$m = 1.67 \times 10^3\text{kg}$,
$I_z = 2.1 \times 10^3\text{kg}\cdot\text{m}^2$,
$l_f = 0.99\text{m},
l_r = 1.7\text{m}$,
$C_{\alpha,f} = 6.1595 \times 10^4\text{N/rad}$, and
$C_{\alpha,r} = 5.2095 \times 10^4\text{N/rad}$.
Fig.~\ref{Fig3} presents the control inputs of the tracking model and the resulting tracking performance. It is clear that the tracking trajectory closely follows the planned trajectory, demonstrating the effectiveness of the learned controller.
To further evaluate the roles of behavior cloning and SMT-based falsification in the proposed framework, we provide additional comparative experiments in Appendix~\ref{apd:fifth}.

\subsection{10D Quadrotor with a 3D Holonomic Planner} 
\label{sec:example:quad}

\begin{figure}[t]
\centering
\includegraphics[width=0.9\textwidth]{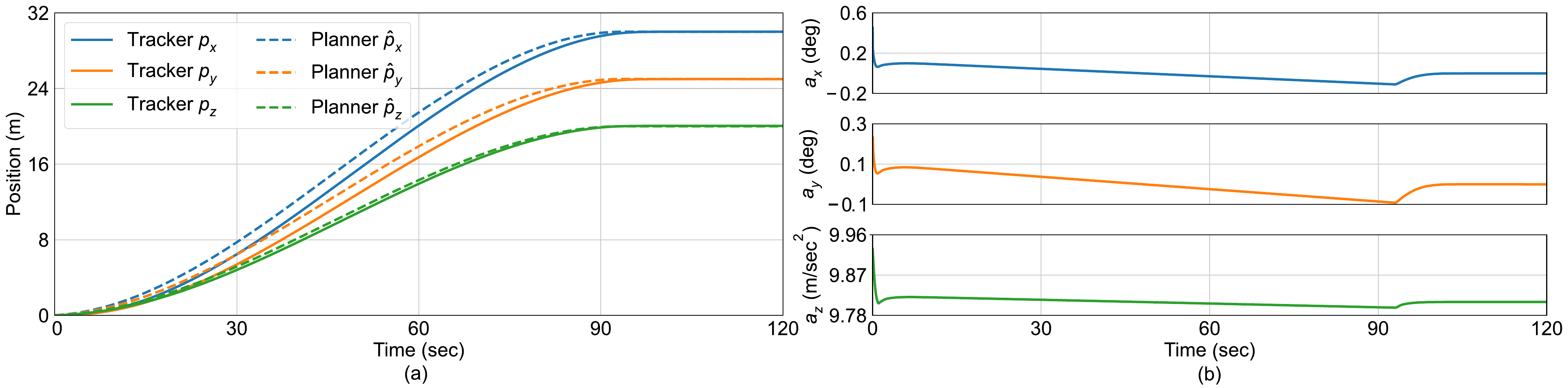}
\caption{(a) Trajectories of the planning and tracking models. (b) Control inputs of the tracking model.}
\label{Fig4}
\end{figure}

\begin{figure}[t]
\centering
\includegraphics[width=0.9\textwidth]{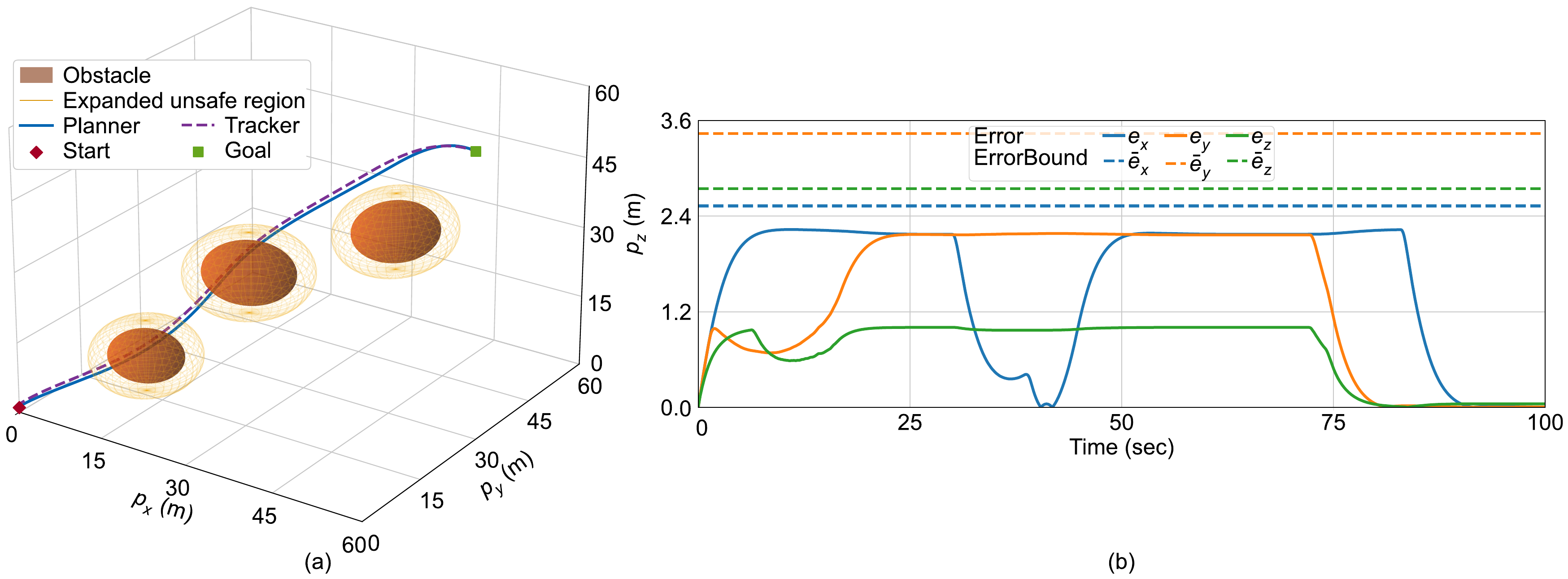}
\caption{(a) Planning and tracking trajectories in the presence of obstacles, where the expanded unsafe regions account for the TEB. (b) Tracking errors and the corresponding error bounds.}
\label{Fig5}
\end{figure}

Consider a $10$-dimensional quadrotor tracking model and a $3$-dimensional holonomic planning model, with their dynamics and constraints given in Appendix~\ref{apd:third} \citep{herbert2017fastrack}.
Fig.~\ref{Fig4} presents the tracking performance and the control inputs of the tracking model in an obstacle-free environment. The results demonstrate  that the learned controller enables the quadrotor to track the planned motion, while satisfying the prescribed input constraints.
Finally, we consider an environment with obstacles to demonstrate the integration of the learned TEB with online planning. As shown in Fig.~\ref{Fig5}(a), the unsafe regions are expanded according to the learned TEB $\mathcal E$, and the planning trajectory generated by MPC avoids the expanded regions. Consequently, the tracking trajectory remains collision-free despite the tracking error. Fig.~\ref{Fig5}(b) further shows that the tracking errors remain within their corresponding bounds.

\section{Conclusion}
\label{sec:conclusion}

This paper presented a learning-based hierarchical planning and tracking framework for nonlinear control-affine systems. We formulated the joint synthesis of a tracking controller, a storage function, and an ellipsoidal TEB as an optimization problem, and developed a learning-based approach for an approximate solution, where the training combined behavior cloning with a counterexample-guided procedure.
For online implementation, the resulting TEB was incorporated into the planner to enforce safety constraints, while the learned controller was applied to the tracking model.
Future work will investigate learning the planning model when it is not known a priori, and validating the proposed framework on real-world robotic systems.


\bibliography{references}

\clearpage
\appendix

\section{Illustrative Examples for Dynamics \eqref{dyn:error}}
\label{apd:first}


\subsection{Vehicle Platoon}

Consider a platoon consisting of a leader vehicle and a follower vehicle. The state of the high-fidelity model is
$x = [x_l, v_l, x_f, v_f]^\top$, where $(x_l, v_l)$ and $(x_f, v_f)$ denote the position and velocity of the leader and follower vehicles.
The corresponding dynamics follows the model (25)-(26) in \cite{smith2018hierarchical}.
To facilitate online planning, we introduce a reduced-order model
$z = [\hat x, \hat v]^\top$, which only describes the motion of the leader vehicle.
Suppose that the desired spacing between the two vehicles is a constant $L > 0$. Then a natural map $\pi$ is
$\pi(z) = [\hat x, \hat v, \hat x - L, \hat v]^\top$.
Furthermore, the error $e: = x - \pi(z)$ measures the deviation of the actual platoon state from the desired configuration.

\subsection{Bicycle Model with a Dubins Planner}
\label{append:A:2}

Following \cite{rajamani2006vehicle}, we consider a tracking system described by the bicycle model
\begin{equation*}
\left\{
\begin{aligned}
\dot X &= v_x \cos\psi - v_y \sin\psi, \\
\dot Y &= v_x \sin\psi + v_y \cos\psi, \\
\dot \psi &= \omega, \\
\dot v_x &= \omega v_y + a_x, \\
\dot v_y &= -\omega v_x + 2/m \cdot \big(F_{c,f}\cos\delta_f + F_{c,r}\big), \\
\dot \omega &= 2/I_z \cdot\big(l_f F_{c,f} - l_r F_{c,r}\big),
\end{aligned}
\right.
\end{equation*}
with
\begin{equation*}
\begin{aligned}
F_{c,f} &= -C_{\alpha,f}\alpha_f,~~
F_{c,r} = -C_{\alpha,r}\alpha_r, \\
\alpha_f &= (v_y + l_f \omega)/v_x -\theta,~~
\alpha_r = (v_y - l_r \omega)/v_x,
\end{aligned}
\end{equation*}
where $X$ and $Y$ are the positions, $\psi$ is the inertial heading angle, $v_x$ and $v_y$ are the longitudinal and lateral speeds, and $\omega$ is the yaw rate.
$m$ and $I_z$ denote the vehicle’s mass and yaw inertia, $l_f$ and $l_r$ represent the distance from the center of mass of the vehicle to the front and rear axles, and $C_{\alpha,i}$ ($i \in \{f,r\}$) is the tire cornering stiffness.
The inputs $\theta$ and $a_x$ are the front wheel steering angle and longitudinal acceleration.

The planner is a Dubin's vehicle model given by
\begin{equation*}
\dot {\hat X} = \hat v \cos\hat\psi,
~\dot {\hat Y}= \hat v \sin\hat\psi, 
~\dot {\hat \psi} = \hat\omega,
\end{equation*}
where $\hat X$ and $\hat Y$ are the positions, and  $\hat \psi$ is the heading angle. The inputs $\hat v$ and $\omega$ denote the speed and angular velocity.

Let $x = [X, Y, \psi, v_x, v_y, \omega]^\top$, and $z = [\hat X, \hat Y, \hat \psi]^\top$. 
For the error variable \eqref{error}, we take $\pi(z) = [z, 0_3]^\top$, and $\phi(x,z) = \text{diag}(R^{-1}(\hat \psi), I_4)$, where
\begin{equation*}
R(\hat \psi) =
\begin{bmatrix} 
\cos(\hat\psi),  -\sin(\hat\psi) \\
\sin(\hat\psi), ~\cos(\hat\psi)
\end{bmatrix}.
\end{equation*}
Note that $\phi$ allows us to replace the trigonometric functions in $\psi$ and $\hat \psi$ by trigonometric functions in $\psi - \hat\psi$.
As a result, the error dynamics can be represented in the form of \eqref{dyn:error}.

\section{Proof of Proposition \ref{prop:Lyap}}\label{apd:second}

\hspace*{1.2em} (i) We prove (a) by contradiction. Suppose that there is $t_2 > 0$ such that $V(e(t_2)) > \gamma$. 
Since $V(e(0)) \le \gamma$, there exists $t_1 \in [0, t_2)$ such that $V(e(t_1)) = \gamma$ and $\frac {d}{dt} V(e(t))|_{t = t_1} \ge 0$ by the continuity of $V$.
This contradicts \eqref{lyap:cond}, and hence, part (a) holds.

(ii) If $e(0) \notin \Omega_\gamma^V$, it follows from \eqref{lyap:cond} that $V(e)$ is nonincreasing whenever $e(t)\notin\Omega_\gamma^V$. 
Moreover, if the trajectory enters $\Omega_\gamma^V$, it remains in $\Omega_\gamma^V$ thereafter by part (a).
Therefore, $V(e(t))$ is bounded for any $e(0) \in \mathbb R^{\tilde n}$. Furthermore, the trajectory $e(t)$ is bounded for all $t \ge 0$ since $V$ is radially unbounded.

By part (a), if there exists $T_1 \ge 0$ such that $e(T_1) \in \Omega_\gamma^V$, then $\text{dist}(e(t), \Omega_\gamma^V) = 0, \forall t \ge T_1$. Therefore, it suffices to consider the case where $e(t) \notin \Omega_\gamma^V$ for all $t \ge 0$.
In this case, $V(e(t))$ is monotonically decreasing and bounded below by $\gamma$. Consequently,
$\ell := \lim_{t\to\infty}V(e(t))$ exists and satisfies $\ell \ge \gamma$.

We claim that $\ell = \gamma$.
Suppose, by contradiction, that $\ell > \gamma$, i.e.,
$\ell-\varepsilon > \gamma$ for some $\varepsilon > 0$.
Since $\lim_{t\to\infty}V(e(t)) = \ell$, there exists $T_2 > 0$ such that $V(e(t))\ge \ell-\varepsilon, \forall t\ge T_2$.
Furthermore, it holds that
$$e(t) \in \widehat \Omega := \{e \in \mathbb R^{\tilde n}: \ell-\varepsilon \le V(e) \le V(e(0))\}, \forall t \ge T_2.$$

It is clear that $\widehat \Omega$ is a compact set. Since $\mathcal V$ and $\mathcal D$ are compact and $\dot V$ is continuous and strictly negative, there exists $c > 0$ such that $\dot V \le -c$, for all $t \ge T_2$, which implies that $V(e(t))$ cannot converge to the finite limit $\ell$.
Thus, $\ell = \gamma$, and moreover, $\lim_{t\to\infty}\operatorname{dist}\bigl(e(t),\Omega_\gamma^V\bigr)=0$.
This completes the proof.

\section{Tracking and Planning Models in Section~\ref{sec:example:quad}} \label{apd:third}

The high-fidelity model is a 10-dimensional quadrotor system of the form
\begin{equation*}
\begin{alignedat}{3}
\dot{p}_x &= v_x + d_x, \qquad&
\dot{v}_x &= g\tan\theta_x, \qquad&
\dot{\theta}_x &= -d_1\theta_x + \omega_x, \\
\dot{p}_y &= v_y + d_y, &
\dot{v}_y &= g\tan\theta_y, &
\dot{\theta}_y &= -d_1\theta_y + \omega_y, \\
\dot{p}_z &= v_z + d_z, &
\dot{v}_z &= k_T a_z - g, &
\dot{\omega}_x &= -d_0\theta_x + n_0 a_x, \\
&&
\dot{\omega}_y &= -d_0\theta_y + n_0 a_y. &&
\end{alignedat}
\end{equation*}
where $(p_x, p_y, p_z)$ denote the position, $(v_x,v_y,v_z)$ denote the velocity, $(\theta_x,\theta_y)$ denote the pitch and roll angles,
and $(\omega_x,\omega_y)$ denote the pitch and roll rates.
The control inputs are $(a_x,a_y,a_z)$, where $a_x$ and $a_y$ represent the desired pitch and roll angles $(|a_x|,|a_y|\le 10 \text{degree}$,
and $a_z$ represents the vertical thrust $(0\le a_z\le 1.5g)$.
The disturbances acting on the tracking model are
$(d_x,d_y,d_z)$, which represent wind disturbances affecting the position dynamics in each direction $(|d_x|,|d_y|,|d_z|\le 0.1~\mathrm{m/s})$.
The parameters are chosen as
$d_0=10$, $d_1=8$, $n_0=10$, $k_T=0.91$, and $g=9.81 {\rm m/s^2}$.

The planning model is given by
\begin{equation*}
\dot{\hat{x}} = b_x,~~
\dot{\hat{y}} = b_y,~~
\dot{\hat{z}} = b_z.
\end{equation*}
The control inputs are $(b_x,b_y,b_z)$, which correspond to the velocities along the three spatial coordinates $(|b_x|,|b_y|,|b_z|\le 0.5\mathrm{m/s})$.

\section{Tracking Performance for the 4D car} \label{apd:forth}

\begin{figure}[t]
\centering
\includegraphics[width=0.9\textwidth]{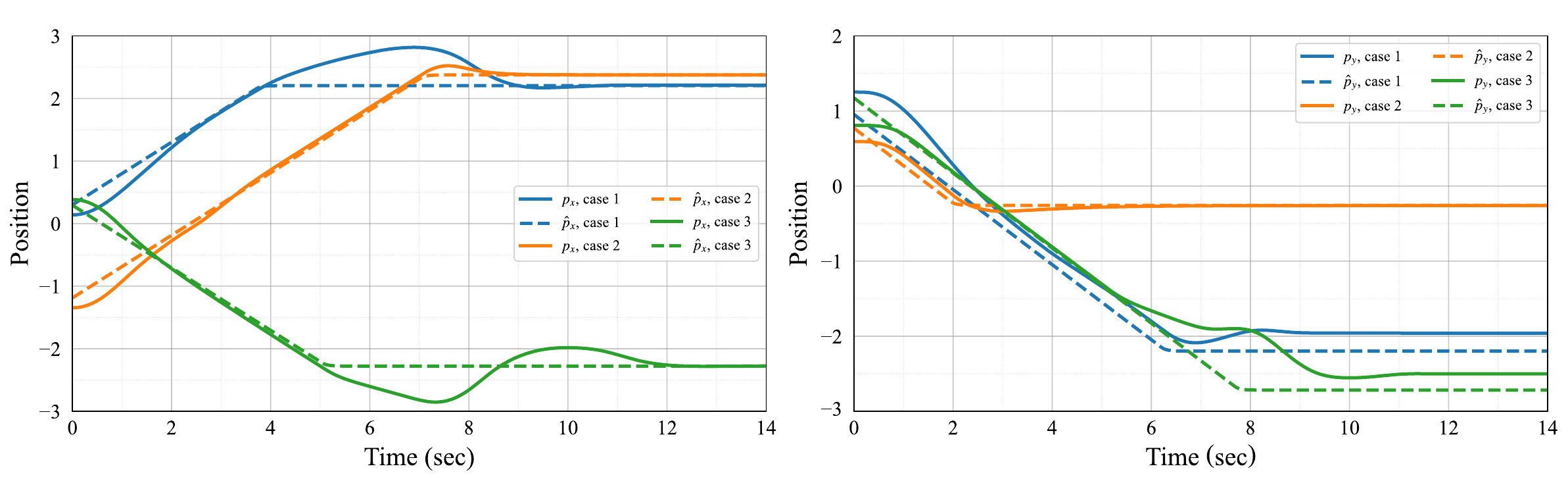}
\caption{Tracking performance of the 4D car model under different goal states.}
\label{Fig6}
\end{figure}

To further evaluate the tracking performance of the learned controller, we consider three additional tasks with different goal states, where the tracking and planning models are initialized at different states. Fig. \ref{Fig6} presents the trajectories of the planning and tracking models along the $p_x$ and $p_y$ coordinates. The results show that the learned controller achieves bounded tracking errors across different planning tasks.

\section{Tracking Performance for the Bicycle Model} \label{apd:fifth}

\begin{figure}[t]
\centering
\includegraphics[width=0.9\textwidth]{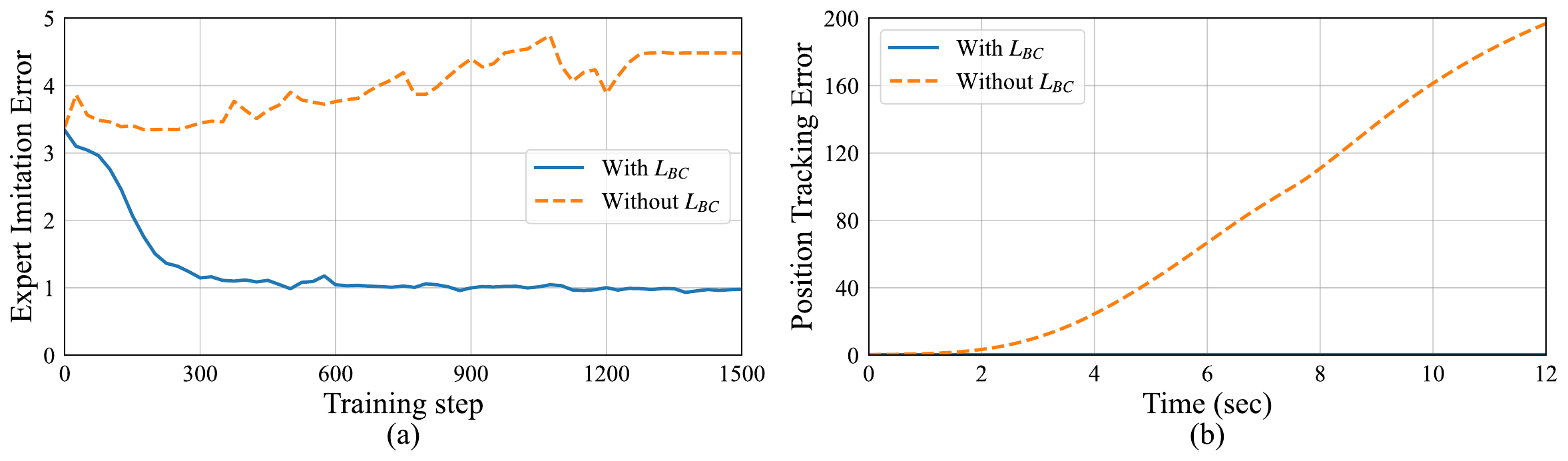}
\caption{(a) Comparison of control input errors relative to the expert policy during training. (b) Comparison of position tracking errors.}
\label{Fig7}
\end{figure}

To evaluate the effect of behavior cloning, we compare the proposed approach with a variant trained without the behavior cloning loss $\mathcal L_{\mathrm{BC}}$ in \eqref{loss}, while keeping all other settings unchanged.
As shown in Fig. \ref{Fig7}(a), incorporating behavior cloning reduces the imitation error.
Fig.~\ref{Fig7}(b) further demonstrates that the controller trained with behavior cloning achieves smaller tracking errors under the same planning task.

\begin{figure}[t]
\centering
\includegraphics[width=0.9\textwidth]{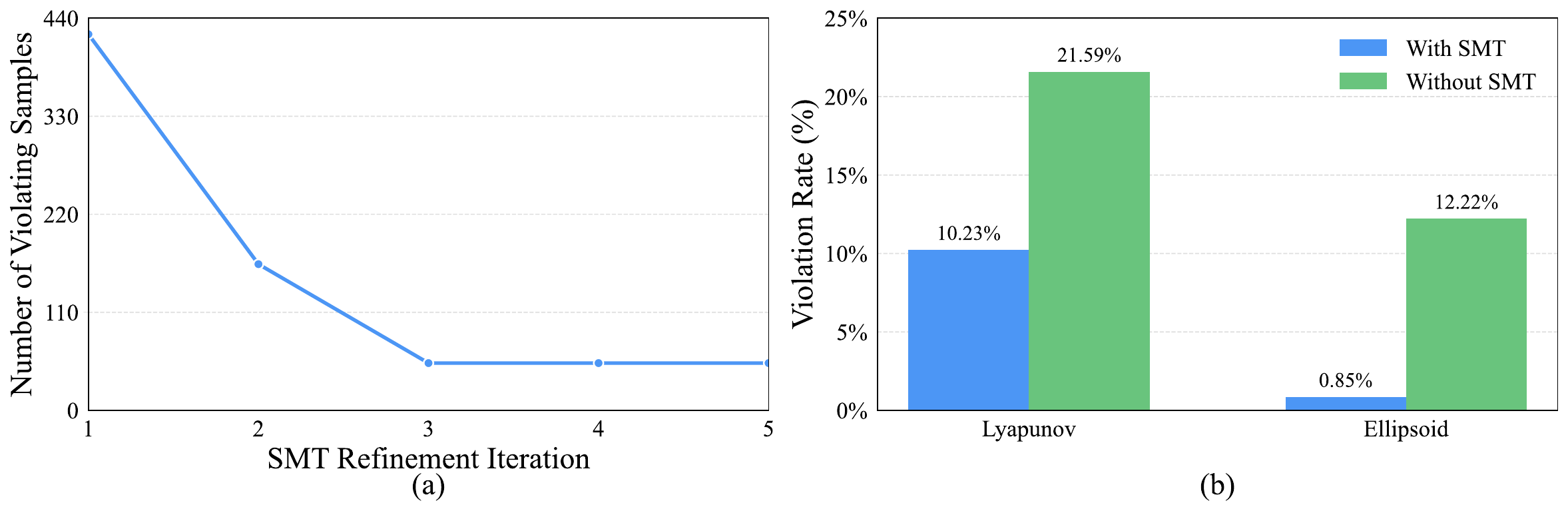}
\caption{(a) Number of samples violating the Lyapunov decrease or ellipsoidal containment conditions during SMT refinement. (b) Comparison of Lyapunov decrease and ellipsoidal containment violation rates.}
\label{Fig8}
\end{figure}

To assess the role of SMT-based falsification, we compare the proposed approach with its counterpart trained without counterexample-guided refinement. Fig.~\ref{Fig8}(a) presents the number of samples that violate the Lyapunov decrease or ellipsoidal containment conditions across SMT refinement iterations.
We further generate an independent evaluation dataset by randomly sampling state and input variables within their prescribed domains. Fig.~\ref{Fig8}(b) shows that incorporating SMT substantially reduces the violation rates.

\end{document}